\documentclass[12pt]{article}
\usepackage{amssymb,amsmath,amsthm,authblk}
\usepackage{braket}
\usepackage{graphicx}
\usepackage{mathrsfs}
\usepackage{subfig}
\usepackage{wrapfig}
\usepackage{amsfonts}
\usepackage[english]{babel}
\usepackage{hyperref}
\usepackage{verbatim} 
\usepackage{standalone}
\usepackage{amsmath}
\usepackage{footmisc}
\usepackage{geometry}
\usepackage{color}

\numberwithin{equation}{section}
\usepackage{enumitem}
\usepackage{relsize}
\usepackage{comment}
\usepackage{xcolor}

\input ArtNouvc.fd

\usepackage{aurical}
\usepackage[T1]{fontenc}
\usepackage{pbsi}
\usepackage[T1]{fontenc}
\usepackage{calligra}
\usepackage[T1]{fontenc}

\usepackage{numprint}

\newcommand{\cG}{{\mathcal G}}

\newcommand{\cM}{{\mathcal M}}

\newcommand{\cQ}{{\mathcal Q}}
\newcommand{\cR}{{\mathcal R}}

\newcommand{\tC}{{\mathtt C}}

\newcommand{\e}{\varepsilon}

\newcommand{\R}{\mathbb R}

\newcommand{\T}{\mathbb T}
\newcommand{\C}{\mathbb C}
\newcommand{\csi}{\xi}

\newcommand{\us}{{s}}
\newcommand{\ur}{{2r}}

\usepackage{amssymb,amsmath,amsfonts,mathrsfs,upgreek,textgreek,bbold,relsize} 
\usepackage{mathabx}

\usepackage[nottoc]{tocbibind}

\usepackage{graphicx} 
\usepackage{version}
\usepackage{mathrsfs}
\usepackage{comment}

\renewcommand\comment[1]{{\iffalse #1 \fi}}

\newcommand\sopra[2]{\genfrac{}{}{0pt}{}{#1}{#2}}

\newtheorem{lemma}{Lemma}[section]
\newtheorem{proposition}{Proposition}[section]

\newtheorem{assumption}{Assumptions}[section]

\newtheorem{remark}{Remark}[section]
\newtheorem{notation}{Notation}[section]
\newtheorem{definition}{Definition}[section]

\newcommand{\real}{ {\mathbb R}   }
\newcommand{\torus}{ {\mathbb T}    }
\newcommand{\integer}{ {\mathbb Z}   }

\newcommand{\complex}{ {\mathbb C}   }

\renewcommand{\Im}{\, {\rm Im}\,}
\renewcommand{\Re}{\, {\rm Re}\,}

\renewcommand{\proof}{\nl {\bf Proof}\  \ }
\newcommand{\eproof}{\qed}%{\hfill $\Box$}

\newcommand\beq[1]{ \begin{equation}\label{#1} }
\newcommand{\eeq}{ \end{equation} }
\newcommand{\beqno}{ \[ }
\newcommand{\eeqno}{ \] }
\newcommand\beqa[1]{ \begin{eqnarray} \label{#1}}
\newcommand{\eeqa}{ \end{eqnarray} }
\newcommand{\beqano}{ \begin{eqnarray*} }
\newcommand{\eeqano}{ \end{eqnarray*} }

\newcommand\dfn[1]{ \begin{definition}\label{#1} }
\newcommand\edfn{ \end{definition} }
\newcommand\ass[1]{ \begin{assumption}\label{#1} }
\newcommand\eass{ \end{assumption} }

\newcommand\notat[1]{ \begin{notation} \label{#1} 
 }
\newcommand\enotat{\end{notation}}

\newcommand\rem{\begin{remark} 
\rm 
}
\newcommand\erem{\end{remark} 
}

\newcommand\equ[1]{{\rm (\ref{#1})}}
\newcommand{\nl}{{\smallskip\noindent}}

\newcommand{\qedeq}{\hskip.5truecm
\vrule width 1.7truemm height 3.5truemm depth 0.truemm}
\renewcommand\qed{\qedeq}

\newcommand{\x}{\xi}

\renewcommand{\a }{\alpha }

\newcommand{\g }{\gamma}

\renewcommand{\l }{\lambda }
\renewcommand{\L }{\Lambda }

\renewcommand{\o }{\omega }

\newcommand{\Z}{\mathbb{Z}}

\renewcommand\AA{{\rm A}}
\newcommand\hAA{\hat\AA}

\newcommand\aaa{{a}}
\newcommand\bbb{{b}}

\newcommand\ttx{{\mathtt x}}
\newcommand\tty{{\mathtt y}}
\newcommand\ttz{{\mathtt z}}
\newcommand\ttu{{\mathtt u}}

\newcommand\ttY{{\mathtt Y}}

\newcommand\ttq{{\mathtt q}}
\newcommand\ttp{{\mathtt p}}

\def\R{\mathbb R}
\def\T{\mathbb T}

\def\bks{\, \backslash\, }
\DeclareMathOperator{\meas}{meas}

\newcommand\eqby[1]{\stackrel{\equ{#1}}{=}}

\newcommand\proiezione{\, { \pi}}

\newcommand\chr\ro

\newcommand{\K}{{\mathtt K}}
\newcommand{\KO}{{\mathtt K}_{\rm o}}

\newcommand{\ro}{{\mathtt r}}

\newcommand\Fio{\Phi_{\!{}_0}}
\newcommand\ttZ{{\tt Z}^k}
\newcommand\ttB{{\tt B}^k}
\newcommand\ttr{{\tt r}}

\newcommand\smin{\underline{s}}
\newcommand\smax{\bar{s}}

\newcommand{\ts}{\textstyle}

\usepackage{appendix}

\def\genKO{{\cal G}^n_{\KO}}
\def\genK{{\cal G}^n_{\K}}

\newcommand\noruno[1]{  |#1|_{{}_1} }

\renewcommand\ln{\log}

\newcommand\pk{{\proiezione^\perp_k}}

\newcommand\Rd{{\mathcal  R}^2}

\newcommand\ddi\updelta			
\newcommand\tk{{\tau_k}}
\newcommand\tktilde{{\tilde{\tau}_k}}

\newcommand\cdr{{c}_{{}_{\star}}}

\renewcommand{\subset}{\subseteq}

\makeatletter
\newcommand{\pushright}[1]{\ifmeasuring@#1\else\omit\hfill$\displaystyle#1$\fi\ignorespaces}
\newcommand{\pushleft}[1]{\ifmeasuring@#1\else\omit$\displaystyle#1$\hfill\fi\ignorespaces}
\makeatother

\begin{document}

\title{{\bf Singular KAM theory for  convex Hamiltonian systems\footnote{L. Biasco and L. Chierchia were supported by 
the grant {\sl NRR-M4C2-I1.1-PRIN 2022-PE1-Stability in Hamiltonian dynamics and beyond-F53D23002730006-Financed by E.U.--NextGenerationEU. }}
\footnote{
Santiago Barbieri was supported by the Juan de la Cierva Postdoctoral Grant JDC2023-052632-I funded by the Spanish National Agency for Research (Agencia Estatal de Investigaci\'on)
}}}

%\date{March 2024}

\author[1]{S. Barbieri}
\author[2]{L. Biasco}
\author[2]{L. Chierchia}
\author[3]{D. Zaccaria}
\affil[1]{\small Universitat Polit\`ecnica de Catalunya, Barcelona (Spain)}
\affil[2]{\small Universit\`a degli Studi Roma Tre, Roma (Italy)}
\affil[3]{\small University of Toronto, Toronto, Ontario (Canada).}

\date{July 12, 2025}

\maketitle

\begin{abstract}\noindent
%REV
In this note, we briefly discuss how the singular KAM Theory of \cite{singularKAM} -- which was worked out for the mechanical case $\frac12 |y|^2+\e f(x)$ --  can be extended  to {\sl convex} real analytic nearly integra\-ble Hamiltonian systems
with Hamiltonian in action--angle variables given by $h(y)+\e f(x)$ with $h$ convex and generic $f$. 

\vskip0.5truecm


\vskip0.1truecm


\end{abstract}

%The extension of Singular KAM theory to generic analytic perturbations (depending also on the action variables) of convex integrable systems is a natural future development of the present work. However, such a task is highly non-trivial, as the constructions in \cite{NL,Complex Arnold,GP,Singular KAM} need to be deeply modified by making use of the quantitative Morse-Sard theory developed by Yomdin and Comte in \cite{Yomdin}. A first step in this direction can be found in \cite{TesiSanti}.

%\tableofcontents

%\section{Introduction}

\nl

\nl
In the Springer {\sl Encyclopaedia of Mathematical Sciences}, Arnold, Kozlov and Neishtadt conjectured that the relative measure of phase space free of invariant tori of a real--analytic, nearly integrable general Hamiltonian System with three or more degrees of freedom is of the {\sl same order of the perturbation}; compare \cite{AKNDUE}. 

\nl
In \cite{singularKAM}, in the special case of natural systems with a Hamiltonian function given (in action--angle variables) by $H=\frac12 |y|^2 + \e f(x)$, it has been proven a related result, namely, that, for generic real--analytic $f$, the relative measure of invariant primary and secondary tori is at least $O(1-\e|\log \e|^c)$ for some $c>0$, in agreement  (up to the logarithmic correction) with Arnold, Kozlov and Neishtadt conjecture\footnote{For related partial results, see \cite{MNT}.}.

\nl
Proofs in \cite{singularKAM} are based on a new ``singular KAM theory'', which extends the classical  theory of Kolmogorov, Arnold and Moser \cite{Ko54,A63,Mo62},  so as to deal, in particular, with primary and secondary tori arbitrarily close to action--angle singularities arising near simple resonances, 
%REV
{\sl uniformly in phase space}.  Let us recall that classical KAM theory predicts that the relative measure of {\sl primary} tory is of order\footnote{Compare \cite{Laz,Nei,P82,Sva}.} $O(1-\sqrt\e)$, and that such an estimate is optimal, as simple integrable cases show.

\nl
In this paper we briefly discuss how to extend  \cite{singularKAM} to {\sl convex} integrable Hamiltonian in place of the purely quadratic term $\frac12 |y|^2$. 
%REV

\nl
The main issue here is geometric: one has to quantitatively describe suitable neighborhoods of simple resonances, so as to be able to put, after averaging, the
secular Hamiltonians in  ``generic standard form''  (as defined in  \cite{GP}). More specifically, we introduce action--angle variables so that the secular Hamiltonian, after averaging over the $(n-1)$ fast angles, depends only on one slow angle. Then, we show that the level sets of the derivative of the integrable part with respect to the slow action, are graphs over {\sl uniform} domains of slow actions; see \S~1.2 below.  At this point, the arguments to show that the density of all KAM primary and secondary tori is at least of order $O(1-\e|\log \e|^c)$ for some $c>0$, follows, without extra difficulties,  the arguments spelled  out in \cite{singularKAM}, which lead to the following statement, confirming (up to the logarithmic correction) the Arnold, Kozolv, Neishdatd conjecture, namely:

\nl
{\sl  For generic\footnote{The generic class we refer to is introduced in \cite{GP}, Sec.1.}  real--analytic potentials $f$ on $\T^n$, the relative measure of invariant primary and secondary tori for 
the Hamiltonian $H=h(y) + \e f(x)$, with $h$ strictly convex, 
is at least $O(1-\e|\log \e|^c)$ for some $c>0$.}

\nl
We also mention a technical improvement in the averaging theory discussed here, which may be useful in applications (for example, in celestial mechanics), namely,  we allow for different analyticity radii in the angle variables; indeed,  in singular KAM theory it is essential to have sharp control of  analytic singularities in the angular variables and having a common width of analyticity strip is quite  unnatural\footnote{For example, in the circular restricted three body problem the angles have different analyticity domains; compare \cite{FZ}.}.

\nl
A second extension would be to replace the perturbation $f$ with a generic function depending also on actions, but this is a much more difficult problem\footnote{For example, it is not obvious how to generalize the generic real--analytic class introduced in \cite{singularKAM} in view of the zeros of $f$ introduced, in general, by the action dependence.
In this respect, as suggested to us by Laurent Niederman, it might be useful to consider the quantitative Morse--Sard theory developed by Yomdin and Comte in \cite{Yomdin}; compare \cite{TesiSanti}.}.

\section{Resonance analysis}

Consider a bounded phase space  $\cM:=B \times \T^n$, with $B \subset \R^n$ a  {\sl bounded convex open
non empty set}, and  $\T^n$  the standard flat $n$--torus  $\R^n / (2 \pi \Z^n)$, endowed with the standard symplectic two--form $dy \wedge dx =\sum_{i=1}^n dy_i \wedge dx_i$. 

\nl
Let $H(y,x)$ be a real--analytic, nearly--integrable  Hamiltonian on $\cM$ given by 
\begin{equation}
    \label{Hamiltonian2}
    H(y,x) = H_{\e}(y,x) =h(y) + \e f(x); \qquad (y,x) \in \cM,\qquad 0\le \e \le 1.
\end{equation} 

\subsubsection*{Analyticity parameters}

We introduce quantitative analyticity parameters as follows.
Let $r>0$ and $\us=(\us_1,...,\us_n) \in \R^n_{+}$ be a vector with positive components. Denote by\footnote{As usual, `bar' denotes complex--conjugated and `dot' inner product.}  $|z|=\sqrt{z\cdot \bar z}$  the standard Euclidean norm on vectors $z\in\complex^n$, and define the following complex neighborhoods 
\beqano
&&    B_{\rho} := \bigcup_{y \in B} \{z \in \C^n : |z-y| <\rho  \};  \\
&&    \T_{\us}^n :=\{x \in \C^n : |\mbox{Im} \ x_i| < \us_i  \ \forall \ i=1,...,n \} / (2 \pi \Z^n).
\eeqano
Henceforth,  {\sl we assume that $H$ is real--analytic and bounded on $B_{2r}\times  \T_{\us}^n$} for some $r>0$.

\rem
As mentioned above, the reason for allowing different analyticity widths in the angular variables is motivated by physical examples, such as the three--body problem, where this is the case. Having sharp control on the angle complex singularities is essential in analyzing fine properties of Hamiltonian systems, especially, in the contex of singular KAM Theory or, possibly, of Arnold diffusion.
\erem

\subsubsection*{Convexity assumption}

In this paper {\sl we  assume that $h$ in \equ{Hamiltonian2} is  $\g$--convex} for some positive $\g$, i.e., we assume that 
    \begin{equation} \label{h-convex}
         \partial^2_y h(y)  \ \csi \cdot \csi = \sum_{i,j=1}^n \partial_{y_i y_j} h(y)\, \csi_i \csi_j \ge \gamma |\csi|^2\,,\qquad\forall \x\in\R^n\,,\ \forall y\in \Re(B_\ur)\,.
    \end{equation}
We shall also assume, without loss of generality, that {\sl  the frequency map $y\to \omega(y)= \partial_y h(y)$ is bi--Lipschitz} and satisfies for some positive constants $\bar L$, $L$ the inequalities
\beq{LIP}
\bar{L}^{-1}  \le\frac{|\o(y)-\o(y_0)|}{|y-y_0|} \le L\,, \ \forall \ y\neq y_0 \in B_\ur\,,\qquad \big(\o(y):= \omega(y)= \partial_y h(y)\big)\,;
\eeq
We also let $M := \sup_{y \in B_\ur} |\omega(y)| < +\infty$.

\subsection{Geometry of resonances and coverings of $B$}
  
In singular KAM Theory \cite{singularKAM}, as well as in Nekhoroshev's Theory (see, e.g., \cite{AKNDUE}), it is fundamental to consider carefully the {\sl geometry of resonances}, which allows for high order averaging theory. 

\nl
We recall, that, given an integer vector $k \in \integer^n$, a \textit{resonance} $\cR_k$ for the  integrable Hamiltonian $h(y)$ is the set $\{ y \in B: \omega(y) \cdot k=0\}$;  a \textit{double resonance} $\cR_{k,\ell}$ is a set given by  $\cR_{k,\ell} = \cR_k \cap \cR_\ell$ with $k$ and $\ell$   linearly independent integer vectors.

\nl
In these definitions, the integer vectors may be assumed to have no common divisors; indeed, it is enough to consider integer vectors
in the set $\cG^n$ defined as {\sl the subset of integer vectors $k \in \Z^n$ with co--prime components and with first non--null component strictly positive\footnote{This set coincides with the set of generators of maximal one--dimensional lattices in $\integer^n$.}}; we shall also denote, for $K>0$,  by 
$\cG^n_K$ the vectors  $k\in \cG^n$ with 1--norm $\noruno{k}:=\sum|k_j|$ not exceeding $K$.

\nl
As well known, not all resonances have to be taken into account, and, typically,
one introduces a ``Fourier cut--off'' $\K$ corresponding to a prefixed ``small--divisor threshold'' $\a>0$, and considers resonances corresponding to Fourier modes of order less or equal than $\K$, the higher ones being negligible in view of the exponentially fast decay of Fourier modes. However, in singular KAM Theory a {\sl double} scale of Fourier modes has to be taken into account, as explained, e.g., in the introduction of \cite{NL}. In the present case the definition of these thresholds have to be slightly modified (with respect to the natural systems case, considered in \cite{singularKAM}) in view of convexity and of the more general assumption on angular analyticity. We therefore give the following definitions. 

\nl
Let $s\in\R^n_+$ be as above, and let $k\in\integer$; we denote
\begin{equation} \label{esse}
    \smin := \min\limits_{1\le i\le n} s_i,
\qquad \smax := \max\limits_{1\le i\le n} s_i\,,\qquad  
\hat s:=\smax/\smin\,,\qquad
|k|_s = \sum_{i=1}^n s_i |k_i| \,.
\end{equation}
Then, following \cite{singularKAM} (compare, Eq. (20), p. 12), we introduce 
two  Fourier scales $\K,\KO$ (to be eventually defined as suitable functions of $\e$) and a  small divisor threshold $\a$ satisfying\footnote{Here, the only (trivial) difference with \cite{singularKAM} is the introduction of the constant $\hat s\ge 1$.}
    \begin{equation}
    \label{K's 2}
    \K \ge 6\hat s\, \KO \ge 6 \KO \ge 12\,,\qquad \alpha:= \sqrt{\e} \K^\nu\,,\quad  \nu := \textstyle \frac{9}{2}n + 2\,.
    \end{equation}
   
 \nl
 Next,  we   define a covering of the action--space, according to a non--resonant zone, simply--resonant zones and  a doubly--resonant zone. Indeed, 
 denote by $\proiezione^\perp_k$
the orthogonal
projection  on the subspace perpendicular  to\footnote{Explicitly, $\proiezione^\perp_k \o:= \o- \frac{1}{|k|^2}(\o\cdot k) k$.} 
$k$,  and
let\footnote{Recall \equ{h-convex} and \equ{LIP}; observe that $L\ge\l_{\rm max}$ and $\g\le \l_{\rm min}$, where  $\l_{\rm max}$ and $\l_{\rm min}$ 
are, respectively, the maximal and minimum eigenvalue of the Hessian of $h$.}
 $$\tC = \tC(n,L,\gamma):=
12\, c_1\,n\,L/\gamma\,,\quad{\rm where}\quad c_1:=5n(n-1)^{n-1}\,;$$
then we define the following subsets of $B$:
\beqa{sonno}
&&\textstyle 
\cR^{0}:= \left\{ y \in B : |\omega(y) \cdot k| \ge \frac{\alpha}{2 \tC} \,, \forall\ \ 0<|k|_{{}_1}\leq \KO\,
\right\}\,,\nonumber\\
&&\textstyle 
\cR^{1k}:= \left\{y \in B: |\o (y)\cdot k|\leq \frac{\a}\tC \,, \mbox{and} \  |\pk \o (y)\cdot \ell|\geq \frac{3 \a \K^{n+3}}{|k|}  \, , \forall
\ell\in \genK\bks\Z k\right\}\,, \nonumber\\
&&
\cR^1 := \bigcup_{k \in \cG^n_{1,\KO}} \cR^{1k}\,, \\
&&
\cR^2 := B \setminus \big( \cR^1 \cup \cR^0 \big)\,.\nonumber
\eeqa

\rem
These definitions are adapted from \cite{singularKAM}: compare, in particular, Eq.'s (21), (22) and (23) at p.~12 of \cite{singularKAM}. 
Notice, first, that in the mechanical case the frequency map is simply the identity map. The second difference is the appearance here of the constant $\tC\ge 1$: this is technical and will become clear below. Finally, in the lower bound on $|\pk \o (y)\cdot \ell|$ in the definition of $\cR^{1k}$ there appears   $\K^{n+3}$ in place of $\K$ (compare  (22) in \cite{singularKAM}); the choice of the power $n+3$ here is not optimal but it is done  for simplicity (as it allows for a single covering of the simply--resonant region), and 
does not affect in any substantial way the strategy of \cite{singularKAM} (which is robust with respect to powers of $\K$, which, at the end, are chosen as  suitable powers of $|\log \e|$).
\erem

\nl
Clearly, from the definition of $\cR^2$ it follows that $\{\cR^i\}$ is a cover of $B$, i.e., that 

$$B= \cR^0\cup\cR^1\cup \cR^2\,.
$$

\nl
The set $\cR^2$ contains all double (or higher) resonances, as well as resonances with high frequency modes. But, as remarked in \cite{singularKAM}, $\cR^2$ is a small set (as we shall shortly see), which is out of reach of perturbation theories\footnote{Compare the Introduction in \cite{singularKAM}.}. 

\nl
Indeed, we claim that the measure of $\cR^2$ can be bounded as
\begin{equation}\label{teheran44} 
\meas (\Rd) \le \cdr\,  \a^2\   \K^{2n+2} \le \cdr \, \e \, \K^{b}\,, 
\end{equation}
where $b=11n + 6$ and $\cdr = 
\frac{c}{2 \cdot 3^n} M^{n-2}\,\bar{L}^n$
for a suitable $c>0$ depending only on $n$.
This is essentially Lemma 2.1 of \cite{singularKAM}; for completeness we reproduce the simple geometrical proof in the convex case.

\nl
\proof {\bf of \equ{teheran44}} Let $\Omega:=\o(B)$ and 
for  $k\in\genKO\,,\ \ell\in \genK\bks\Z k$, 
define 
\beqa{defi}\textstyle
&&\textstyle \Omega^2_{k,\ell} := \left\{ \o\in\Omega: |\omega \cdot k| < \frac{\alpha}{\tC}   \ {\rm and}\  |\pk \omega \cdot \ell| \le \frac{3 \alpha \K^{n+3} }{|k|}\right\}\,;
\nonumber\\
&&\Rd_{k\ell}:= \left\{y\in B: \o \in \Omega^2_{k,\ell} \right\}\,.
\eeqa
Then, one checks easily that 
\beq{due}
\Rd \subseteq\bigcup_{k\in \genKO} 
\bigcup_{ \sopra{\ell\in \mathcal G^n_{\K}}{\ell\notin  \Z k}}  \Rd_{k,\ell}\,.
\eeq
Now, 
denote by  $v\in\R^n$  the projection of $\o$ onto  
the plane generated by $k$ and $\ell$
(recall that, by hypothesis, $k$ and $\ell$ are not parallel);
then,   
\begin{equation}\label{soldatino}\textstyle
|v\cdot k|=|\o\cdot k|<\a\,, \qquad |\proiezione_k^\perp v \cdot \ell|
=|\proiezione_k^\perp \o \cdot \ell|
 \le 
\frac{3\a\K^{n+3}}{|k|}\,.
\end{equation}
Set
%REV. h -> \bar\ell
\beq{bacca}\ts
\bar \ell:=\pk \ell= \ell -\frac{\ell\cdot k}{|k|^2} k\,.
\eeq
Then, $v$ decomposes in a unique way as
$v=a k+ b \bar \ell$
for suitable $a,b\in\R$.
By \eqref{soldatino},
\beq{goja}\textstyle
|a|<\frac{\a}{|k|^2}\,,\qquad
|\pk v\cdot\ell|
=|b\bar\ell\cdot \ell| \le \frac{3\a\K^{n+3}}{|k|}\,,
\eeq
and, since $ |\ell|^2 |k|^2-(\ell\cdot k)^2$ is a positive integer (recall, that $k$ and $\ell$ are  not parallel), 
$$
|\bar\ell\cdot \ell|
\eqby{bacca}
\frac{ |\ell|^2 |k|^2-(\ell\cdot k)^2  }{|k|^2}
\ge \frac1{|k|^2}\,.
$$
Hence, 
\beq{velazquez}\textstyle
|b|\le \frac{3 \a \K^{n+3} }{|k|} \,.
\eeq
Then,  write $\o \in \Omega^2_{k,\ell}$ as $\o=v+v^\perp$ with 
$v^\perp$ in the orthogonal
complement of the plane generated by $k$ and $\ell$. Since $|v^\perp |\le |\o|< M$  and $v$ lies in the plane spanned by $k$ and $\ell$ inside a rectangle of sizes of length $2\a/|k|^2$ and $6 \a \K^{n+3} |k|$ (compare \equ{goja} and \equ{velazquez}),
we find that, for any $k\in\genKO$ and $\ell\in  \mathcal G^n_{\K}\bks \Z k$, one has 
\beqno\ts
\meas(\Omega^2_{k,\ell})\le \frac{2\a}{|k|^2}\, (6 \a \K^{n+3} |k|)\ (2M)^{n-2}=3\cdot 2^n \, M^{n-2} \, \a^2 \frac{\K^{n+3}}{|k|}\,.
\eeqno
Since $\sum_{k\in\genKO}|k|^{-1}\le c\, \KO^{n-1}$ for a suitable $c=c(n)$,  by \equ{LIP}, 
Eq.~\ref{teheran44} follows. \qed

\subsection{Simple resonances as graphs over adiabatic actions}

As already mentioned, a key r\^ole in singular KAM theory is played by simple resonances. It is therefore important to have a precise analytic description of them. In this subsection, we express neighborhoods of simple resonances as a foliation of graphs over $(n-1)$ adiabatic actions; see Proposition~\ref{BWV208} below. This is the first step to put the averaged secular Hamiltonian in the so--called ``generic standard form''; compare \S 3 in \cite{GP}.

\nl
To state the results in Proposition~\ref{BWV208}, we need some preparation.

\nl
Let $k\in\cG^n$.
A simple resonances $\o\cdot k=0$ corresponds to a  ``resonant angle'',  in the sense that one  can make a linear symplectic change of variables
on $B\times \torus^n$ so that in the new variables one transforms the resonant relation $k\cdot x$ into, say, the first new angle  $\ttx_1=k\cdot x$. Indeed, one can find  a matrix $\AA\in {\rm SL}(n,\Z)$ (i.e., integer--valued with determinant one) with first row given by   $k$; such a matrix $\AA$ can be chosen so that\footnote{Compare Lemma 2.6, (i) in \cite{GP}; $|M|_{{}_\infty}$, with $M$ matrix (or vector), denotes the maximum norm $\max_{ij}|M_{ij}|$ (or $\max_i |M_i|$). }
\beq{atlantide}
 \AA:=\binom{k}{\hAA}\,,\ 
|\hAA|_{{}_\infty}\leq |k|_{{}_\infty}\,,\ \ 
|\AA|_{{}_\infty}=|k|_{{}_\infty}\,,\ \ 
|\AA^{-1}|_{{}_\infty}\leq 
(n-1)^{\frac{n-1}2} |k|_{{}_\infty}^{n-1}\,. 
\eeq
Then, the linear symplectic change of variables
\begin{equation}\label{fiodena}
\Fio: (\tty,\ttx) \mapsto (y,x)= (\AA^T\tty,  \AA^{-1} \ttx)\,,\qquad \quad (\ttx_1=k\cdot x)\,.
\end{equation} 
transforms the unperturbed Hamiltonian $h(y)$ into\footnote{Notice that also ${\rm B}$ is convex.}
\beq{dsa}
\mathtt h^k_0 (\tty):= h( \AA^T \tty)\,,\qquad \tty\in \ttB:=\AA^{-T}B\,.
\eeq
Note that $\mathtt h^k_0$ is holomorphic 
on the complex set\footnote{By
\eqref{atlantide} $\|A\|:=\max_{|v|=1}|Av|
\leq n|A|_\infty=n|k|_\infty$.} 
\begin{equation}\label{HWV432}
\ttB_{2\ttr_k}\,,\quad\mbox{with}\quad  \ttr_k:=
\frac{r}{n|k|_\infty}\,.
\end{equation}
%and the following estimates follow from  \eqref{LIP} and\footnote{By
%\eqref{atlantide} 
%$\|A^{-1}\|\leq n(n-1)^{\frac{n-1}2} |k|_{{}_\infty}^{n-1}$; then $|Av|\geq |v|/\|A^{-1}\|\geq |v|/n(n-1)^{\frac{n-1}2} |k|_{{}_\infty}^{n-1}$.}
% \eqref{atlantide}    \begin{equation}\att{
%    \begin{split}\label{LIP2} 
%    &\big(\bar{L}c_2  |k|_{{}_\infty}^{n-1}\big)^{-1} |y-y_0| \le|\nabla
%    \mathtt h^k_0(y)
%    -\nabla
%    \mathtt h^k_0(y_0)| 
%    \le n|k|_\infty L |y-y_0|, \qquad \forall \ \tty,\tty_0 \in \ttB_{2\ttr_k},
%    \\
%    &\sup_{\ttB_{2\ttr_k}} |\nabla
%    \mathtt h^k_0 | \leq n|k|_\infty M\,,
%    \end{split}
%    }
%    \end{equation}
%where $c_2=n(n-1)^{\frac{n-1}2}$. 

%REV
\begin{notation}\label{vettori}
We shall adopt the following conventions: a vector $\tty\in\real^n$ will be denoted
$$
\tty=(\tty_1,\tty_2,...,\tty_n)=(\tty_1,\hat \tty)\,,\qquad \hat \tty=(\tty_2,....,\tty_n)\in\real^{n-1}\,,
$$
namely, the hat over $n$--vectors denote the projection onto the coordinates with index greater or equal than two.
\end{notation}

\noindent
Now, since
\begin{equation}\label{ventaglietti}
\partial_{\tty_1} \mathtt h^k_0 (\tty)=
\partial_yh(\AA^T \tty)\cdot k=
\o(\AA^T \tty)\cdot k
\end{equation}
 and noting that
the first column of $\AA^T$ is exactly $k$ (by \eqref{atlantide}), by \eqref{h-convex}
and \eqref{LIP} we obtain
\beq{alessandro}
|\partial_{\tty_1} \mathtt h^k_0 (\tty)
    -\partial_{\tty_1} \mathtt h^k_0 (\tty_0)| 
    \le 
    L |k|  |\tty-\tty_0|, \qquad \forall \ \tty,\tty_0 \in \ttB_{2\ttr_k}
 \eeq
 and, by convexity, for  $(\tty_1,\hat\tty),(\tty_1',\hat\tty) \in \ttB_{2\ttr_k}$, 
 \beq{magno}
 \partial_{\tty_1} \mathtt h^k_0 (\tty_1,\hat\tty)
    -\partial_{\tty_1} \mathtt h^k_0 
    (\tty_1',\hat\tty) \geq
    \g |k|^2 (\tty_1-\tty_1')
    , \quad \forall \ 
    \tty_1>\tty_1'\,.
\eeq
Indeed, by Lagrange's theorem, there exists a point
$\ttz_1$
between $\tty_1$ and $\tty_1'$ such that\footnote{Note that $(\ttz_1,\hat\tty) \in \ttB_{2\ttr_k}$ since
$\ttB_{2\ttr_k}$ is a convex set.}
\begin{eqnarray*}
   \partial_{\tty_1} \mathtt h^k_0 (\tty_1,\hat\tty)
    -\partial_{\tty_1} \mathtt h^k_0 
    (\tty_1',\hat\tty)
    &=&
    \partial^2_{\tty_1} \mathtt h^k_0 (\ttz_1,\hat\tty)
     (\tty_1-\tty_1')\\
     &\stackrel{\eqref{dsa}}=&
   \Big(\partial^2_y  h\big( \AA^T (\ttz_1,\hat\tty)\big)k\cdot k\Big)
   (\tty_1-\tty_1')\nonumber\\
  & \stackrel{\eqref{h-convex}}\geq& \g |k|^2(\tty_1-\tty_1').
\end{eqnarray*}
In these new variables
 the simple resonance $\o(y)\cdot k=0$ becomes $\partial_{\tty_1} \mathtt h^k_0 (\tty)=0$. 
 
 \nl
 We therefore define a (suitable) real neighborhood of the simple
resonance in $\ttB$ by letting
\beq{balletto}
\ttZ:=\{\tty\in\Re(\ttB_{\frac54\ttr_k}):
 \partial_{\tty_1} \mathtt h^k_0 (\tty)=0\}\,.
\eeq
We also set 
\beq{balletto2}
\ttZ_\varpi:=\{\tty\in\Re(\ttB_{\ttr_k}):
 \partial_{\tty_1} \mathtt h^k_0 (\tty)=\varpi\}\,,
\eeq
so that $\ttZ\supseteq\ttZ_0$.

\nl
Now, for a fixed  $k\in \cG^n$,  one may express (because of convexity)  the resonant hyper--surface  $\ttZ$ in 
\eqref{balletto} as a graph over the last
$n-1$ actions $\hat \tty:=(\tty_2,\ldots,\tty_n)$, finding $\tty_1=\eta(\hat \tty)$.

\nl
Let 
\begin{equation}\label{rkgotico}
\mathfrak r_k:=
\frac{\tktilde r}{8n^{3/2}|k|}
\leq 
\frac{\ttr_k}{8\sqrt n},
 \quad
 \tk:=\frac{\g|k|}{2L},
 \quad
\tktilde:=\min\{1,\tk\}.
\end{equation}
and define the following  $(n-1)$--dimensional cubes of edge $\mathfrak r_k$ and lower left corner in 
%REV
$\mathfrak r_k j$ with $j\in\integer^{n-1}$:
\beq{Qj}
\cQ_j:=\mathfrak r_k\cdot \big(j+[0,1)^{n-1}\big)
\subset \mathbb R^{n-1}\,,
\qquad
j\in\mathbb Z^{n-1}.
\eeq
Note that $\cQ_j\cap \cQ_{j'}=\emptyset$
if $j\neq j'$, while $\bigsqcup_{j\in\mathbb Z^{n-1}}\cQ_j=\mathbb R^{n-1}$.
Define the projection of a set 
$S\subset \mathbb R^n$
on a set $\hat E\subset R^{n-1}$ as
$$
\Pi_{\hat E}S:=\{\hat \tty\in \hat E \ :\ 
\exists \tty_1\in\mathbb R\ \mbox{with}\ 
(\tty_1,\hat\tty)\in S\}.
$$
Set
\begin{equation}\label{def:J}
J:=\{
j\in\mathbb Z^{n-1} \ :\ 
\Pi_{\cQ_j}\ttZ\neq\emptyset
\}.
\end{equation}
Since $\ttZ$ is bounded,
$J$ is a finite set.
Note that
\begin{equation}\label{parco}
\Pi_{\mathbb R^{n-1}}\ttZ
\subset
\cQ:=
\bigsqcup_{j\in J}\cQ_j
\subset
\Pi_{\mathbb R^{n-1}}
\big(\Re(\ttB_{\frac32\ttr_k})\big).
\end{equation}

\begin{proposition}\label{BWV208}
 There exists a real analytic function
 \begin{equation}\label{felicita}
 \eta_0^k:[-\varpi_0^k,\varpi_0^k]\times\cQ\subset\R\times \R^{n-1}\to\R\,,
 \qquad
 \varpi_0^k:=\sqrt n L|k| \mathfrak r_k,
\end{equation} 
 such that
 \begin{equation}\label{cacioepepe}
  \partial_{\tty_1} \mathtt h^k_0 \big(
 \eta_0^k(\varpi,\hat\tty),\hat\tty\big)=\varpi
 \end{equation}
 and, for any $\varpi\in[-\varpi_0^k,\varpi_0^k]$,
\begin{equation}\label{peles}
 \ttZ_\varpi\subseteq 
 \tilde{\tt Z}_\varpi^k:=\big\{\big(\eta_0^k(\varpi,\hat\tty),\hat\tty\big)\ :\ 
 \hat\tty\in\cQ\big\}
  \subset
 \Re(\ttB_{\frac32\ttr_k}).
\end{equation}
 Moreover 
 \begin{equation}\label{mare}
\sup_{[-\varpi_0^k,\varpi_0^k]\times\cQ}
|\partial_\varpi \eta_0^k|
\leq \frac{1}{\g |k|^2}.
\end{equation}
 Finally $\eta_0^k(0,\cdot)$
 possesses 
 holomorphic extension 
 on\footnote{Obviously
 one could holomorphically extend $\eta_0^k$
 also in its first variable, but we need here 
 only the extension in the second variable.}  $\cQ_{\hat\ttr_k}$,
 with, for $0<\ttr\leq \ttr_k$,
\begin{equation}\label{diana}
 \sup_{\cQ_{\hat\ttr_k}}|\Im\eta_0^k(0,\cdot)|
 \leq \frac{1}{\tk}\hat\ttr_k,
 \qquad
 \hat\ttr_k:=
 \frac{1}{2^9 n}
\tktilde^2\ttr.
 \end{equation}
 
\end{proposition}
\proof
We first construct the function $\eta^k_0$.
Fix $\varpi\in[-\varpi_0^k,\varpi_0^k]$ and 
$\hat\tty\in \cQ_j$ for some $j\in J$ (recall \equ{Qj}).
By \eqref{balletto} and the definition of $J$ in \eqref{def:J}
we know that 
 there exists a given
$$\ttz=(\ttz_1,\hat\ttz)\in \Re(\ttB_{5\ttr_k/4})$$
with  $\hat\ttz\in\cQ_j$ and 
$\partial_{\tty_1} \mathtt h^k_0 (\ttz)=0$.
Since $\hat\tty,\hat\ttz\in\cQ_j$, 
by \eqref{rkgotico} we have that
\begin{equation}\label{nonna}
(\ttz_1,\hat\tty)\in \Re(\ttB_{11\ttr_k/8})\,.
\end{equation}
 Setting
\begin{equation}\label{rkgoticotilde}
 \tilde{\mathfrak r}_k:=
 2\frac{\sqrt n L}{\g |k|}\mathfrak r_k
 \stackrel{\eqref{rkgotico}}\leq
 \frac{r}{8n|k|}
 \stackrel{\eqref{HWV432}}\leq 
 \frac{\ttr_k}{8},
 \end{equation}
 we have that 
  the segment
$
[\ttz_1-\tilde{\mathfrak r}_k,
 \ttz_1+\tilde{\mathfrak r}_k]\times\{\hat\tty\}
$
belongs to the convex set $\Re(\ttB_{\frac32\ttr_k})$,
namely
\begin{equation}\label{sweetdreams}
[\ttz_1-\tilde{\mathfrak r}_k,
 \ttz_1+\tilde{\mathfrak r}_k]\times\{\hat\tty\}
 \subset
 \Re(\ttB_{\frac32\ttr_k}).
\end{equation}

Since,
 by  \eqref{alessandro},
we get
 $$
 | \partial_{\tty_1} \mathtt h^k_0 (\ttz_1,\hat\tty)|
 =
 | \partial_{\tty_1} \mathtt h^k_0 (\ttz_1,\hat\tty)-
 \partial_{\tty_1} \mathtt h^k_0 (\ttz_1,\hat\ttz)|
 \leq
 L|k||\hat\tty-\hat\ttz|
 \leq
 \sqrt n L|k| \mathfrak r_k
  \,,
 $$
we have that
$$
\partial_{\tty_1} \mathtt h^k_0 
(\ttz_1+\tilde{\mathfrak r}_k,\hat\tty)
-\varpi
\stackrel{\eqref{magno}}\geq
-\varpi_0^k+
\g |k|^2 \tilde{\mathfrak r}_k
+
\partial_{\tty_1} \mathtt h^k_0 (\ttz_1,\hat\tty)
\geq 
-\varpi_0^k+
\g |k|^2 \tilde{\mathfrak r}_k
- \sqrt n L|k| \mathfrak r_k
= 0,
$$
by \eqref{felicita} and \eqref{rkgoticotilde}.
Analogously 
$\partial_{\tty_1} \mathtt h^k_0 
(\ttz_1-\tilde{\mathfrak r}_k,\hat\tty)-\varpi\leq 0$.
Note that, since the domain $\Re(\ttB_{2\ttr_k})$
is convex, the function 
$\tty_1\to \partial_{\tty_1} \mathtt h^k_0 
(\tty_1,\hat\tty)$
is defined on a interval, which contains $[\ttz_1-\tilde{\mathfrak r}_k,
 \ttz_1+\tilde{\mathfrak r}_k]$; moreover it is 
 continuous and strictly increasing
by \eqref{magno}.
As a consequence, 
 there exists a unique value
$$\aaa\in[\ttz_1-\tilde{\mathfrak r}_k,
 \ttz_1+\tilde{\mathfrak r}_k]$$ such that
 $\partial_{\tty_1} \mathtt h^k_0 
(\aaa,\hat\tty)=\varpi$.
Then we set
$\eta_0^k(\varpi,\hat\tty):=\aaa$.
Recollecting we have defined 
$\eta_0^k$ in \eqref{felicita}, satisfying
\eqref{cacioepepe}.
Moreover the last inclusion in \eqref{peles}
follows by \eqref{sweetdreams}.
The fact that $\eta^k_0$ is real analytic
follows since $\mathtt h^k_0$ is
real analytic.

\nl
We now prove the first inclusion in \eqref{peles}.
Fix $\varpi\in[-\varpi_0^k,\varpi_0^k]$.
Recalling \eqref{balletto2}, let us take a point
$$\ttu=(\ttu_1,\hat\ttu)\in
\ttZ_\varpi\,,
$$
 namely $\ttu\in\Re(\ttB_{\ttr_k})$ and
 $\partial_{\tty_1} \mathtt h^k_0 (\ttu)=\varpi$.
 If $\hat\ttu\in \cQ$, by unicity $\ttu_1=
 \eta_0^k(\varpi,\hat\ttu)$,
 and, therefore, $\ttu\in  \tilde{\tt Z}_\varpi^k$.
On the other hand it is not possible that $\hat\ttu\notin \cQ$.
Indeed, assume, e.g., that $\varpi=0$; then
$\ttu\in \ttZ_0\subseteq\ttZ$  and, by \eqref{parco},
$\hat\ttu\in \cQ$.
On the other hand, when $\varpi$ is negative, 
since, by \eqref{rkgoticotilde},
 we have that 
  the segment
$$
[\ttu_1,
 \ttu_1+\tilde{\mathfrak r}_k]\times\{\hat\ttu\}
$$
belongs to the convex set $\Re(\ttB_{\frac54\ttr_k})$.
Moreover
$$
\partial_{\tty_1} \mathtt h^k_0 
(\ttu_1+\tilde{\mathfrak r}_k,\hat\ttu)
\stackrel{\eqref{magno}}\geq
\g |k|^2 \tilde{\mathfrak r}_k
+
\partial_{\tty_1} \mathtt h^k_0 (\ttu_1,\hat\ttu)
=\g |k|^2 \tilde{\mathfrak r}_k+\varpi
\geq 
\g |k|^2 \tilde{\mathfrak r}_k-\varpi_0^k
> 0,
$$
by \eqref{felicita} and \eqref{rkgoticotilde}.
Again by continuity there exists 
a value $\bbb$ such that
$$
\partial_{\tty_1} \mathtt h^k_0 (\bbb,\hat\ttu)=0
$$
Since $(\bbb,\hat\ttu)\in\Re(\ttB_{\frac54\ttr_k})$
we have that $(\bbb,\hat\ttu)\in \ttZ$,
then, by \eqref{parco}, $\hat\ttu\in \cQ$.
The case when $\varpi$ is positive is analogous.
We conclude that $\hat\ttu\in \cQ$ in all cases.
The proof of the inclusion \eqref{peles} is completed.

\nl
We finally show that 
$\eta_0^k(0,\cdot)$ has
 holomorphic extension on  $\cQ_{\hat\ttr_k}$
 and satisfies the estimate in \eqref{diana}.
 Fix a point $\tty^0=(\tty_1^0,\hat\tty^0)=
 \big(\eta_0^k(0,\hat\tty^0),\hat\tty^0\big)\in 
 \tilde{\tt Z}^k_0$ with $\hat\tty^0\in\cQ\subset
\Pi_{\mathbb R^{n-1}}
\big(\Re(\ttB_{\frac32\ttr_k})\big)$ by \eqref{parco}.
 By construction 
 $\partial_{\tty_1} \mathtt h^k_0(\tty_1^0,\hat\tty^0)=0$.
 Let $\ttY_1$ be the complex closed ball of radius
\begin{equation}\label{r1}
 \ttr_1:= \frac{1}{\tk}\hat\ttr_k,
 \end{equation}
  centred at $\tty_1^0$
 and let
 $\hat\ttY$ be the complex closed ball of radius
 $ \hat\ttr_k$  centred at $\hat\tty^0$.
 Note that by \eqref{diana}
 \begin{equation}\label{tortellino}
  \ttr_1,\hat\ttr_k
  \leq
   \frac{1}{2^9 n}
 \tktilde\ttr
 \leq
 \frac{1}{2^9 n}
\tktilde\ttr_k
 \leq
 \frac{1}{2^9 n}
\ttr_k.
 \end{equation}
 Since, by \eqref{peles}, $\tty^0\in 
 \tilde{\tt Z}^k_0\subset
 \Re(\ttB_{\frac32\ttr_k})$, we have
 \begin{equation}\label{peperoni}
\ttY_1\times \hat\ttY
\subset \ttB_{\frac74\ttr_k}.
\end{equation}
 Let $E$ be the Banach space of the 
 continuous functions $\eta:\hat\ttY\to \C$
 that are holomorphic in the interior of 
$\hat\ttY$, endowed with the sup-norm.
Let $\mathcal C$ be its closed subset formed by the functions
$\eta:\hat\ttY\to \ttY_1$.
We claim that the map
$$
\eta(\cdot)\to \eta(\cdot)- 
\partial_{\tty_1} \mathtt h^k_0(\eta(\cdot),\cdot)/d\,,\qquad
d:=\partial^2_{\tty_1} \mathtt h^k_0(\tty^0)
$$
is a contraction on $\mathcal C$.
The fixed point of the above map is the 
required (local) holomorphic  extension of $\eta_0^k(0,\cdot)$,
proving \eqref{diana} by the definition 
of  $\ttr_1$.
Since, by \eqref{magno}, $d\geq \g |k|^2$,
it is immediate to see  that the fact that the above map 
is a contraction on $\mathcal C$ follows from the following two estimates:
\begin{equation}\label{scarlatti}
\sup_{\hat\tty\in\hat\ttY}
|\partial_{\tty_1} \mathtt h^k_0(\tty_1^0,\hat\tty)|
\leq \frac12 \g |k|^2\ttr_1\,,
\qquad
\sup_{\tty_1\in\ttY_1,\hat\tty\in\hat\ttY}
|\partial_{\tty_1}^2 \mathtt h^k_0(\tty)-
\partial_{\tty_1}^2 \mathtt h^k_0(\tty^0)|
\leq \frac12 \g |k|^2.
\end{equation}
Note also that 
 the estimate in \eqref{diana}        
 follows  since the image of $\eta_0^k$
 is contained, by construction, in $\ttY_1$
 and by the definition of $\ttr_1$. 
 
 \nl
 It remains to prove         \eqref{scarlatti}.  
  Since $\partial_{\tty_1} \mathtt h^k_0(\tty^0)=0$,
  for any $\tty\in\ttB_{2\ttr_k}$ we obtain
  \begin{equation}\label{raviolo}
 |\partial_{\tty_1} \mathtt h^k_0(\tty)|
 =
 |\partial_{\tty_1} \mathtt h^k_0(\tty)
 -\partial_{\tty_1} \mathtt h^k_0(\tty^0)|
 \stackrel{\eqref{alessandro}}
 \leq
 L |k|  |\tty-\tty_0|,
 \end{equation}
 in particular
 \begin{equation}\label{raviolo2}
 \sup_{\ttB_{2\ttr_k}}|\partial_{\tty_1} \mathtt h^k_0(\tty)|
 \leq
 2L |k| \ttr_k.
 \end{equation}
 Then  the first estimate in \eqref{scarlatti}
 follows by \eqref{raviolo} since, for
 $\hat\tty\in\hat\ttY$,
 $$
 |\partial_{\tty_1} \mathtt h^k_0(\tty_1^0,\hat\tty)|
 \leq
 L |k|  |\hat\tty-\hat\tty_0|
 \leq
 L |k| \hat\ttr_k
\stackrel{\eqref{diana},\eqref{r1}}= \frac12 \g |k|^2\ttr_1.
 $$
 Let us finally prove the second estimate in \eqref{scarlatti}.
 Fix
 $\tty=(\tty_1,\hat\tty)$ with
  $\tty_1\in\ttY_1$ and $\hat\tty\in\hat\ttY$; 
 by \eqref{raviolo2}, \eqref{peperoni}
  and Cauchy estimates we obtain
 \begin{eqnarray*}
|\partial_{\tty_1}^2 \mathtt h^k_0(\tty)-
\partial_{\tty_1}^2 \mathtt h^k_0(\tty^0)|
\leq 
64 n L |k|\frac{\ttr_1+\hat\ttr_k}{\ttr_k}
\stackrel{\eqref{tortellino}}\leq
 \frac{1}{8}L|k|\left(\frac{1}{\tk}+1\right)
\tktilde^2
\stackrel{\eqref{diana}}\leq
\frac12 \g |k|^2,
\end{eqnarray*}
concluding the proof of the second estimate
in \eqref{scarlatti}.

Since 
$$
\partial_\varpi \eta_0^k(\varpi,\hat\tty)=
\frac{1}{\partial_{\tty_1}^2 \mathtt h^k_0
\big(\eta_0^k(\varpi,\hat\tty),\hat\tty\big)}
$$
and, by \eqref{magno} 
$\partial_{\tty_1}^2 \mathtt h^k_0
\geq
\g |k|^2$, we get \eqref{mare}.
\qed											
\bigskip

\nl
Recalling \eqref{sonno} and \eqref{dsa}
we have that
\beqa{tornimparte}
D^k&:=&A^{-T}\cR^{1k}\\
&=&\textstyle  \big\{\tty \in \ttB:
 |\o (A^{T}\tty)\cdot k|\leq \frac\a\tC \,, \mbox{and} \  |\pk \o (A^{T}\tty)\cdot \ell|\geq \frac{3 \a \K^{n+3}}{|k|}  \, , \forall
\ell\in \genK\bks\Z k\big\}\nonumber
\eeqa
and, by using \eqref{ventaglietti},
 \eqref{cacioepepe} and \eqref{peles}, we obtain
\beqa{tornimparte2}
D^k&=& \textstyle
\big\{\tty \in \ttB:
\eta_0^k(-\frac\a\tC,\hat\tty)
\leq
\tty_1
\leq
\eta_0^k(\frac\a\tC,\hat\tty)
 \,,\\
&&\textstyle
\phantom{AAAAA}  \mbox{and} \  |\pk \o (A^{T}\tty)\cdot \ell|\geq \frac{3 \a \K^{n+3}}{|k|}  \, , \forall
\ell\in \genK\bks\Z k\big\}.\nonumber 
\eeqa
 Let us define the normal set
 \begin{equation}\label{assergi}
 \textstyle
\check D^k:=\big\{\tty=(\tty_1,\hat\tty)\ :\ 
\eta_0^k(-\frac\a\tC,\hat\tty)
\leq
\tty_1
\leq
\eta_0^k(\frac\a\tC,\hat\tty)\,,\ 
\hat\tty\in \hat D^k\big\}\,,\quad
\hat D^k:=\Pi_{\R^{n-1}}D^k\,.
\end{equation}
It is obvious that 
\begin{equation}\label{roccacalascio}
D^k\subseteq \check D^k.
\end{equation}

\nl
We recall that, given lattice $\L\subset\Z^n$, one says that
a set $D\subseteq \real^n$ is $(\a,K)$ non--resonant modulo $\L$ for a Hamiltonian $h(y)$ defined on $D$, when 
$|\partial_y h(y)\cdot k|\geq \a$
for every $0<|k|\leq K$, $k\notin\L$, $y\in D$.

\begin{lemma}\label{ostia}
For $\K$ large enough, the set  $\check D^k$ is 
$\big(\frac{2 \a \K^{n+3}}{|k|},\K\big)$
non--resonant modulo $\L:=\Z(1,0,\ldots,0)$.
\end{lemma}
\proof
Let $\tty=(\tty_1,\hat\tty)\in \check D^k$.
Then  
$\eta_0^k(-\frac\a\tC,\hat\tty)
\leq
\tty_1
\leq
\eta_0^k(\frac\a\tC,\hat\tty)$
and 
there exists $\ttz_1$ with 
$\eta_0^k(-\frac\a\tC,\hat\tty)
\leq
\ttz_1
\leq
\eta_0^k(\frac\a\tC,\hat\tty)$
such that
$\ttz=(\ttz_1,\hat\tty)\in D^k$.
For any $\ell\in \genK\bks\Z k$
we have that 
$|\pk \o (A^{T}\ttz)\cdot \ell|\geq
 \frac{3 \a \K^{n+3}}{|k|}$.
 Since 
 $$
 |\ttz-\tty|=|\ttz_1-\tty_1|
 \leq 
 \eta_0^k(\frac\a\tC,\hat\tty)
 -
 \eta_0^k(-\frac\a\tC,\hat\tty)
 \stackrel{\eqref{mare}}\leq
 \frac{2\a}{\g|k|^2\tC},
 $$  
 by \eqref{LIP} we get
 $$
 |\pk \o (A^{T}\tty)\cdot \ell|\geq
 \frac{2 \a \K^{n+3}}{|k|},
 $$
for $\K$ large enough.
\eproof
 
\section{Outline of proofs}

We are now ready to outline the main steps needed to extend singular KAM theory, which, as a consequence, yields the announced lower bound on the measure of primary and secondary tori for generic perturbation.  
 
 \nl
{\bf (1)} Applying the symplectic transformation
 $ \Fio$ in 
 \eqref{fiodena}
 to the Hamiltonian $H$ in \eqref{Hamiltonian2}
one obtains (recall \eqref{dsa}), 
\beq{brendel}
 \mathtt H(\tty,\ttx)=
 \mathtt h^k_0(\tty)+\e f(\AA^{-1} \ttx)\,.
 \eeq
Now, one can do high order averaging theory. More precisely, 
by Lemma \ref{ostia}, $\check D^k$ in \eqref{assergi}
is $\big(\frac{2 \a \K^{n+3}}{|k|},\K\big)$
non--resonant modulo $$\L:=\Z(1,0,\ldots,0)\,,$$
and, taking
  $\K$ suitably large with $\e$ small, 
  e.g. $\K\sim |\ln \e|^2$,
  via a  close-to-the-identity
symplectic transformation defined
 in a suitable complex
neighbourhood  of $\check D^k\times \T^n$, 
one can put the Hamiltonian \equ{brendel}
in normal form:
$$
 \mathtt h^k_0(\tty)+\e \big( g^k_{\rm o}(\tty)+ 
 f^k_1(\ttx_1)+
g^k(\tty,\ttx_1) +
f^k (\tty,\ttx)\big),
$$
where:
\begin{itemize}

\item
$
g^k(\tty,\ttx_1)=o(1)\,,$

\item
$f^k_1(\ttx_1)$ is the 
average  
of
$f(\AA^{-1} \ttx)$ with respect to the fast angles
$\ttx_2,\ldots,\ttx_n$
and $o(1)$ going to 0 as $\e\to 0$;

\item
$f^k$ is (almost\footnote{I.e. smaller than
  any power of $\e$.})
  exponentially small.
  
\end{itemize}  
For details, see, e.g.,   section 2 of \cite{GP}. 
 
 \nl
 {\bf (2)}  Disregarding the  exponentially small
 term  $f^k$,  we consider, now,  the ``effective Hamiltonian''
 $$
\mathtt h^k_0(\tty)+\e \big( g^k_{\rm o}(\tty)+ 
f^k_1(\ttx_1)+
g^k(\tty,\ttx_1)\big),
$$ 
which is a one degree--of--freedom  Hamiltonian (in the dynamic variables $\tty_1$ and $\ttx_1$) 
depending on the ``dumb actions'' 
$$\hat\tty:=(\tty_2,...,\tty_n)
$$
as parameters.
Then, one can construct a  symplectic transformation 
 $$
 \tty_1=\ttp_1+\tilde\eta_0^k(\hat\ttp),
 \qquad
 \ttx_1=\ttq_1,
 \qquad 
 \hat\tty=\hat\ttp,
 $$
 for a suitable function\footnote{Rrecall Proposition \ref{BWV208}.} 
 $\tilde\eta_0^k(\hat\ttp)=\eta_0^k(0,\hat\ttp)
 +O(\e)$, such that,
 in the new variables, the effective  Hamiltonian
 takes the form
 $$
\mathrm h^k(\ttp)+\e\big(f^k_1(\ttq_1)+
\mathrm g^k(\ttp,\ttq_1)\big),
 $$
 where\footnote{Recall \eqref{cacioepepe}
 with $\varpi=0$.}
 \begin{equation}\label{anniversario}
 \partial_{\ttp_1}\mathrm h^k(0,\hat\ttp)=0
 \end{equation}
 and
  $\mathrm g^k=o(1)$.
  
  \nl
 Note that the above transformation can be 
 sympectically completed in the variables 
 $\hat\ttx=\hat\ttq-
 \ttq_1 \partial_{\hat\ttp} \tilde\eta_0^k(\hat\ttp)$,
 however this transformation is not well defined
 on $\T^n$, since it is obviously 
 not periodic in $\ttq_1$.
 The way to overcome this problem is explained in section 3 of \cite{singularKAM}.
 
 \nl
{\bf (3)}  Fixing an arbitrary point $\hat\ttp_0\in \hat D^k$
 (compare \equ{assergi})
 and Taylor expanding the Hamiltonian at $\ttp_1=0$, one gets, by \equ{anniversario},
 \begin{equation}\label{bruchi}
 \mathrm h^k(0,\hat\ttp)+
 \mathrm m_k
 \big(1+O(|\hat\ttp-\hat\ttp_0|)+O(|\ttp_1|)\big)
 \ttp^2
 +\e\big(f^k_1(\ttq_1)+
\mathrm g^k(\ttp,\ttq_1)\big),
 \end{equation}
 where, by convexity,
 $$
 \mathrm m_k:=
 \frac12\partial_{\ttp_1}^2\mathrm h^k(0,\hat\ttp_0)>0.
 $$
It is easy to see that,
up to the inessential term
$\mathrm h^k(0,\hat\ttp)$
and rescaling by $\mathrm m_k$,
 the Hamiltonian in \equ{bruchi},
 for $|\hat\ttp-\hat\ttp_0|$ small, 
 can be put into the ``standard form''
 $$
  \big(1+\upnu(\ttp,\ttq_1)
  \big)
 \ttp^2
 +\e\big(G_0(\ttq_1)+
G(\ttp,\ttq_1)\big)
 $$
 where 
 $$
 G_0(\ttq_1):=f^k_1(\ttq_1)/\mathrm m_k
 $$
 and $\upnu$ and $G$ are small;
 compare section 2.2 of \cite{singularKAM}.
This Hamiltonian is suitable for action angle
variables as discussed in \cite{Complex Arnold}.

\nl
At this point no further technical difficulties arise
and the singular KAM Theory applies as in \cite{singularKAM}, leading to the announced measure estimates on primary and secondary tori.

%%%%%%%%%%%%%%%%%%%%%%%%%%%%%%%%%%%%%%%%%%%%%%%%%%%%%%%%%%%%%%%%%%%%%

\end{document}